\documentclass[10pt,a4paper,oneside]{amsart}
\usepackage{amsmath} \usepackage{amssymb}
\usepackage{latexsym} \usepackage{mathrsfs} 
\usepackage{fancyhdr} 
\usepackage{mathrsfs}
\usepackage{amscd}
\usepackage[all]{xy}

\newtheorem{defi}{\bf Definition}[section]

\newtheorem{lemm}[defi]{\bf Lemma}
\newtheorem{theo}[defi]{\bf Theorem}
\newtheorem{prop}[defi]{\bf Proposition}
\newtheorem{rema}[defi]{\bf Remark}

\newtheorem{claim}{\bf Claim}

\newcommand{\qd}[1]{\rightline{#1 $\Box $}}
\newcommand{\wee}{{\scriptscriptstyle \vee}}
\newcommand{\Db}{\mathfrak{Db}}
\newcommand{\Dbrd}[1]{\mathfrak{Db}_{\widetilde{X_t}}^{{\rm rd}, #1}}
\newcommand{\Crdt}{{\mathcal{C}^{\rm rd}_{\widetilde{X_t}}}}
\newcommand{\Crdtp}[1]{{\mathcal{C}^{{\rm rd}, #1}_{\widetilde{X_t}}}}
\newcommand{\widebar}[1]{\overline{#1}} 
\newcommand{\setmin}{\! \smallsetminus \!}
\newcommand{\R}{{\mathbb R}} \newcommand{\C}{{\mathbb C}}
\newcommand{\N}{{\mathbb N}} \newcommand{\Q}{{\mathbb Q}}
\newcommand{\Ct}{\mathcal{C}_{\widetilde{X_t}}} 
\newcommand{\Ctd}{\mathcal{C}_{\widetilde{X_t}, \widetilde{D_t}}}
\newcommand{\Z}{{\mathbb Z}} 
\newcommand{\Hy}{{\mathbb H}}
\newcommand{\Hom}{{\rm Hom}} 
\newcommand{\cHom}{\mathcal{H}om}

\newcommand{\Epq}[2]{{\Omega^{\infty, (#1,#2)}_{X_t}}}
\newcommand{\Etop}[1]{{\Omega^{\infty,#1}_{\widetilde{X_t}}}}
\newcommand{\sD}{\ast D_t}
\newcommand{\vi}{\varphi} 
\newcommand{\bve}{\mathcal{E}} 
\newcommand{\vt}{\vartheta}
\newcommand{\Int}{\int\limits}
\newcommand{\lto}{\longrightarrow} 
\newcommand{\Crd}{\mathcal{C}^{\rm rd}} 
\newcommand{\Hrd}{H^{rd}} 
\newcommand{\Hdr}{H_{dR}} 

 \renewcommand{\H}{\mathcal{H}}
\renewcommand{\O}{\mathcal{O}} 

\newcommand{\claimpf}{\noindent {\em Proof.} --- }
\newcommand{\prf}[1]{{\noindent {\em Proof of #1.} --- }}
\newcommand{\AD}{\mathcal{A}_{\widetilde{X_t}}^{<D_t}} 
\newcommand{\DR}{{\rm DR}_{\Xan}} 
 
\newcommand{\DRD}{{\rm DR}^{<D_t}_{\widetilde{X_t}}} 
\newcommand{\Amod}{\mathcal{A}_{\widetilde{X_t}}^{{\rm mod }D_t}} 
\newcommand{\DRmod}{{\rm DR}_{\widetilde{X_t}}^{{\rm mod }D_t}}
\newcommand{\DRmd}[1]{{\rm DR}_{\widetilde{X_t}}^{{\rm mod }D_t,#1}}

\newcommand{\bR}{\boldsymbol{R}}
\newcommand{\As}[2]{\mathcal{A}_{\widetilde{#1}}^{#2}}
\newcommand{\Ass}[3]{\mathcal{A}_{\widetilde{#1},#2}^{#3}}

\newcommand{\Smod}{\mathcal{S}^{{\rm mod }D_t}}
\newcommand{\SD}{{}^\wee\!{\mathcal{S}}^{<D_t}}
\newcommand{\dolb}{\widebar{\partial}}
\newcommand{\bdot}{\raisebox{-.3ex}{\mbox{\boldmath $\cdot $}}}
\newcommand{\Rd}{\mathcal{R}d}
\newcommand{\Uan}{U_t^{\rm an}}
\newcommand{\Xan}{X_t^{\rm an}}
\newcommand{\an}{{\rm an}}
\newcommand{\OU}{\O_{U_t}}

\newcommand{\cN}{\mathcal{N}}
\newcommand{\cD}{\mathcal{D}}
\newcommand{\cH}{\mathcal{H}}
 
\newcommand{\Af}{{\mathbb A}}

\renewcommand{\labelenumi}{{\rm \roman{enumi})}}

\makeatletter
\renewcommand{\section}{\@startsection{section}{1}{\z@}%
{-\baselineskip}{1\baselineskip}{\large\bf}}
\renewcommand{\subsection}{\@startsection{subsection}{2}{\z@}%
{-\baselineskip}{0.5\baselineskip}{\normalsize\bf}}
\renewcommand{\subsubsection}{\@startsection{subsubsection}{3}{\z@}%
{-\baselineskip}{0.5\baselineskip}{\normalsize\bf}}
\makeatother

\numberwithin{equation}{section} 

\author[Hien]{Marco Hien}
\address{Marco Hien\\ NWF I -- Mathematik\\ Universit{\"a}t Regensburg\\ 93040
  Regensburg\\ Germany}
\email{marco.hien@mathematik.uni-regensburg.de}

\author[Roucairol]{Celine Roucairol}
\address{Celine Roucairol\\ Lehrstuhl f\"ur Mathematik VI\\ Universit\"at Mannheim\\
A5 6\\ 68161 Mannheim\\ Germany}
\email{celine.roucairol@uni-mannheim.de}

\title[Integral representations]{Integral representations for
  solutions of exponential Gau\ss-Manin systems}

\begin{document}



\begin{abstract}
  Let $f,g:U \to {\mathbb A}^1 $ be two algebraically independent
  regular functions from the smooth affine complex variety $U $ to the
  affine line. The associated exponential Gau{\ss}-Manin systems on the
  affine line are defined to be the cohomology sheaves of the direct image of the exponential
  differential system $\mathcal{O}_U e^g $ with respect to $f $. We
  prove that its holomorphic solutions admit representations in terms
  of period integrals over topological chains with possibly closed support and with rapid decay
  condition.
\end{abstract}

\maketitle

\section{Introduction}

Let $U $ be a smooth affine complex variety of dimension $n $ and let
$f,g: U \to \Af^1 $ be two regular functions which we assume to be
algebraically independent. We consider the flat algebraic connection
$\nabla: \O_U \lto \Omega^1_U $ on the trivial line bundle $\O_U $
defined as $\nabla u= du + dg \cdot u $. Let us denote the associated
holonomic $\cD_U $-module by $\O_U e^g $, i.e. the trivial $\O_U
$-module on which any vector field $\xi $ on $U $ acts via the
associated derivation $\nabla_\xi $ induced by the connection. Let
$\cN $ be the direct image $\cN := f_+(\O_U e^g) $ in the theory of
$\cD $-modules (see e.g.~\cite{borel}). Then $\cN $ is a complex of
$\cD_{\Af^1} $-modules and we consider its $k $-th cohomology sheaf
$$
\cH^k f_+(\O_U e^g) \ ,
$$
a holonomic $\cD_{\Af^1} $-module.

There is an alternative point of view of $\cH^k \cN $ in terms of flat
connections: outside a finite subset $\Sigma_1 \subset \Af^1 $, the
$\O_{\Af^1} $-module $\cH^k f_+(\O_U e^g) $ is locally free, hence a
flat connection. This connections coincides with the Gau{\ss}-Manin
connection on the relative de Rham cohomology
$$
\cH^k(\cN)|_{\Af^1 \smallsetminus \Sigma_1} \cong \big( \bR^{k+n-1}
f_\ast (\Omega_{U|\Af^1}^{\bdot}, \nabla ), \nabla_{GM} \big) |_{\Af^1
  \smallsetminus \Sigma_1} \ ,
$$
i.e. the corresponding higher direct image of the complex of relative
differential forms on $U $ with respect to $f $ and the differential
induced by the absolute connection, as in \cite{katz}. It is a flat
connection over $\Af^1 \smallsetminus \Sigma_1 $.

In the case $g=0 $, which is classically called the Gau{\ss}-Manin
system of $f $, it is well-known that the solutions admit integral
representations (see \cite{pham} and \cite{pham2}). The aim of the present article is to give such a
description of the solutions in the more general case of the
Gau{\ss}-Manin system of the exponential module $\O_U e^g $, which we
will call {\em exponential Gau{\ss}-Manin system}. A major difficulty
lies in the definition of the integration involved. In the case of
vanishing $g $, the connection is regular singular at infinity and the
topological cycles over which the integrations are performed can be
chosen with compact support inside the affine variety $U $. In the
exponential case $g \neq 0 $, the connection is irregular singular at
infinity and we will have to consider integration paths approaching
the irregular locus at infinity.

A systematic examination of period integrals in more general irregular
singular situations over complex surfaces has been carried out in
\cite{hien} resulting in a perfect dualtiy between the de Rham
cohomology and some homology groups, the {\em rapid decay homology
  groups}, in terms of period integrals. We will prove in this article
that this result generalizes to arbitrary dimension in case of an
elementary exponential connection $\O e^g $ as before, see Theorem
\ref{thm:dual}.

With this tool at hand, we construct a local system $\widetilde{H^{rd}_k} $ on
$\Af^1 \smallsetminus \Sigma_2 $, the stalk at point $t $ of which is
the rapid decay homology $\Hrd_k(f^{-1}(t), e^{-g_t}) $ of the
restriction of $\O_U e^{-g} $ to the fibre $f^{-1}(t) $.
The set $\Sigma_2$ can be described using a compactification $(F,G):X\to\mathbb{P}^1\times\mathbb{P}^1$ of the map $(f,g):U\to\mathbb{A}^1\times\mathbb{A}^1$ (see Proposition \ref{hrd flat}). Let $\overline{\Delta}$ be the projective variety in $\mathbb{P}^1\times\mathbb{P}^1$ such that $(F,G)$ is a locally trivial fibration out of $\overline{\Delta}$. The set $\Sigma_2$ can be choosen as the set of $t_0\in\mathbb{A}^1$ such that $(t_0,\infty)$ belongs to the closure of $\overline{\Delta}\setminus(\mathbb{P}^1\times\{\infty\})$ in $\mathbb{P}^1\times\mathbb{P}^1$.

Let $\Sigma=\Sigma_1\cup\Sigma_2$. We prove that given a flat section $c_t \otimes e^{g_t} $ in the
local system $H^{rd}_{k+n-1} $ and a relative differential form $\omega $ in
$\Omega_{rel}^{k+n-1}( f^{-1}( \Af^1 \smallsetminus \Sigma)) $ (which
describes the analytification of $\cH^k \cN $, see Proposition
\ref{prop:analytG}), the integral
\begin{equation} \label{eq:integralintro}
\int_{c_t} \omega|_{f^{-1}(t)} \cdot e^{g_t}
\end{equation}
gives a (multivalued) holomorphic solution of the exponential
Gau{\ss}-Manin system. As a consequence of the duality proved in
Theorem \ref{thm:dual}, we deduce that all (multivalued) holomorphic
solutions on $\Af^1 \smallsetminus \Sigma $ can be obtained by this
construction, i.e. we achieve the following theorem (Theorem
\ref{solutions}) which is the main
result of this article:
\begin{theo}
  For any simply connected open subset $V \subset \C \smallsetminus
  \Sigma $, the space of holomorphic solutions of the exponential
  Gau{\ss}-Manin system $\cH^k \cN $ is isomorphic to the space of
  flat sections in the local system $\widetilde{H^{rd}_{k+n-1}} $ over $V $.
\end{theo}
The isomorphism in the theorem is given in terms of the integration
\eqref{eq:integralintro}. We want to remark that this theorem
generalizes similar results by F.~Pham for the Fourier-Laplace transform of the Gau\ss-Manin systems of a regular function $h:U\to\mathbb{A}^1$ (consider $f:\mathbb{A}^1\times U\to\mathbb{A}^1$ to be the first canonical projection and $g:\mathbb{A}^1\times U\to\mathbb{A}^1$ given by $g(x,u)=xh(u)$, cf.~\cite{pham3}). It also
includes other well-known examples as the integral representations of
the Bessel-functions for instance (cf.~the introduction of \cite{b-e}).

\section{The period pairing}

Consider the situation described in the introduction, namely $f,g:U
\to \Af^1 $ being two algebraically independent regular functions on
the smooth affine complex variety $U $. The main result of this work
is to give a representation of the Gau\ss-Manin solutions in terms of
period integrals. The proof relies on a duality statement between the
algebraic de Rham cohomology of the exponential connection associated
to $g $ and some Betti homology groups with decay condition. The
present section is devoted to the proof of this duality statement,
which is a generalization of the analogous result for surfaces in
\cite{hien} to the case of exponential line bundles in arbitrary
dimensions.

\subsection{Definition}

Let us start with any smooth affine complex variety $U_t $ and a regular
function $g_t: U_t \to \Af^1 $. Note, that the index $t $ is
irrelevant in this section but will become meaningful later in the
application to the Gau\ss-Manin system where $U_t $ will denote the
smooth fibres $f^{-1}(t) $. 

We assume that $U_t $ is embedded into a smooth projective variety $X_t $
with the complement $D_t :=X_t \smallsetminus U_t $ being a divisor
with normal crossings. Moreover, we require that $g_t $ extends to a
meromorphic mapping on $X_t $. In particular, if $y \in D_t $ is a
point at infinity, we can choose coordinates $x_1, \ldots, x_n $ of
$X_t $ centered at $y $ such that locally $D_t=\{ x_1 \cdots x_k =0 \}
$ for some $k \in \{ 1, \ldots, n \} $ and such that we have
$$
g_t(x) = x_1^{n_1} \cdots x_r^{n_r} \cdot u(x)
$$
locally at $y $, where $r \le k $ and the exponents $n_i \in \Z $ are
either all positive or all negative and $u $ is holomorphic and
non-vanishing on $D_t $ (i.e. $g_t $ will not have indeterminancies at
$D_t $). 

The function $g_t $ defines a flat meromorphic connection $\nabla_t $
of exponential type on $U_t $, denoted $\O e^{g_t} $ on the trivial
line bundle $\O_{U_t} $ with connection defined as
$$
\nabla_t 1 = dg_t \in \Omega^1_{U_t|\C} \ .
$$

In \cite{hien}, a duality pairing between the de Rham cohomology of a
flat meromorphic connection (admitting a good formal structure) and a
certain homology theory, the rapid decay homology, is constructed in
the case $\dim(X_t)=2 $, generalizing previous constructions of Bloch
and Esnault for curves. For connections of exponential type lying in
the main focus of this work, we will now generalized these results to
the case of arbitrary dimension. We remark that a crucial point in
\cite{hien} is to work with a {\em good} compactification with respect
to the given vector bundle and connection which is always fulfilled
for the connections of exponential type considered above.

We recall the definition of rapid decay homology as in \cite{hien}.
Consider the real oriented blow-up $\pi:\widetilde{X_t} \to X^{\rm an}_t $ of the
irreducible components of $D_t $ (cp. \cite{majima}). Let $\Ct^{-p} $
denote the sheaf associated to the presheaf $V \mapsto
S_p(\widetilde{X_t}, \widetilde{X_t} \smallsetminus V) $ of $\Q $-vector
spaces, where $S_p(Y) $ denotes the groups of piecewise smooth
singular $p $-chains in $Y $ and $S_p(Y,A) $ the relative chains for a
pair of topological spaces $A \subset Y $.

We adopt the standard sign convention for the resulting complex
$\Ct^{-\bdot} $. i.e. the differential of $\Ct^{-\bdot} $ will be
given by $(-1)^r $ times the topological boundary operator $\partial $
on $\Ct^{-r} $. If $d:=\dim_\C(X_t) $, the complex $\Ct^{-\bdot} $ is a
resolution of the sheaf $\widetilde{\jmath}_!  \C_{\Uan}[2d] $ which
in turn is the dualizing sheaf on the compact real manifold
$\widetilde{X_t} $ with boundary $\widetilde{D_t}=\widetilde{X_t} \smallsetminus \Uan $
(cp.~\cite{verdier}).

Now, let $\Ctd^{-\bdot} $ be the complex of chains relative to the
boundary $\widetilde{D_t} $, i.e. the sheaf associated to
$$
V \mapsto S_{\bdot}(\widetilde{X_t}, (\widetilde{X_t} \smallsetminus V)
\cup \widetilde{D_t}) \ .
$$
Let $c \in \Gamma(V, \Ctd^{-\bdot}) $ be a local section over some
open $V \subset \widetilde{X_t} $, represented by a piecewise smooth map
from the standard $p $-simplex $\Delta^p $ to $\widetilde{X_t} $. We
consider $c $ always together with the solution $e^{g_t} $ of the dual of the given
connection of exponential type and will therefore write $c \otimes e^{g_t}
$ for $c $ in the following definition (in the terminology of
\cite{hien}, we consider $c \otimes e^{g_t} $ as a section in
$\Ctd^{-\bdot} \otimes_\Q \widetilde{\jmath}_\ast \bve $ where $\bve $
denoted the local system associated to the dual bundle $\mathcal{O}_{U_t}e^{-g_t}$, in our case the
trivial system with basis $e^{g_t} $).

Now, consider a point $y \in c(\Delta^p) \cap \widetilde{D_t} \cap V $.
Choose local coordinates $x_1, \ldots, x_d $ of $X_t $ around $y=0 $
such that $D_t = \{ x_1 \cdots x_k = 0 \} $.

\begin{defi} \label{def:rdchains} We call $c \otimes e^{g_t} $ a {\em
    rapid decay chain}, if for all $y $ as above, the function
  $e^{g_t(x)} $ has rapid decay for the argument approaching
  $\widetilde{D_t} $, i.e.~if for all $N \in\N^k $ there is a $C_N > 0 $
  such that
  $$
  |e^{g_t(x)}| \le C_N \cdot |x_1|^{N_1} \cdots |x_k|^{N_k}
  $$
  for all $x \in (c(\Delta^p) \setmin \widetilde{D_t}) \cap V $ with
  small $|x_1|, \ldots, |x_k| $.
\end{defi}
In other terms, the requirement is that $\arg(g_t(x)) \in
(\frac{\pi}{2}, \frac{3 \pi}{2}) $ for all $x $ on $c $ with small
$|x| $.

In case $V \cap \widetilde{D_t}=\emptyset $, we do not impose any
condition on $c \otimes e^{g_t} $. Let $\Crdtp{-p}(e^{-g_t}) $ denote the
subsheaf of $\Ctd^{-p} \otimes e^{g_t} $ generated by all rapid decay $p
$-chains and all chains inside $U_t $. As for the sign notation, note
that the dual exponential module is denoted by $\O e^{-g_t} $ and
admits $e^{g_t} $ as holomorphic solution. Together with the usual boundary
operator of chains, these give the {\bf complex of rapid decay chains}
$\Crdt(e^{-g_t}):=(\Crdtp{-\bdot}(e^{-g_t}), \partial) $.

\begin{defi}
  The {\bf rapid decay homology} of $\O e^{-g_t} $ is the hypercohomology
  $$
  \Hrd_k(U_t, e^{-g_t}) := \Hy^{-k}(\widetilde{X_t}, \Crdt(e^{-g_t})) \ ,
  $$
  where $\pi:\widetilde{X_t} \to \Xan $ denotes the oriented real
  blow-up of the normal crossing divisor $D_t:= X_t \smallsetminus U_t
  $.
\end{defi}
In the same way as in \cite{hien}, Prop. 3.10, one proves that this
definition does not depend on the choice of the compactification $X_t $.
We will comment on this later.

The usual barycentric subdivision operator on $\Ctd $ induces a
subdivision operator on the rapid decay complex. The latter is
therefore easily seen to be a homotopically fine complex of sheaves
(cp.~\cite{swan}, p.~87). Hence, the rapid decay homology can be
computed as the cohomology of the cooresponding complex of global
sections:
$$
\Hrd_k(U_t, e^{-g_t}) = H^{-k} \left( \Gamma(\widetilde{X_t},
  \Crdt(e^{-g_t})) \right) \ .
$$

In other words, we can regard the elements of this group as relative topological cycles $c$ in $H^k(X_t,D_t,\mathbb{C})$ which satisfy the rapid decay condition for $g_t$.

Since $U_t $ is assumed to be affine, the de Rham cohomology of $\OU
e^{g_t} $ can also be computed on the level of global sections
$$
\Hdr^p(U_t, e^{g_t}) := \Hy^p( U_t, (\Omega_{U_t|k}^{\bdot}, \nabla_t))
\cong H^p(
\ldots \to \Gamma_{U_t}(\Omega_{U_t|k}^q) \stackrel{\nabla_t}{\to}
\Gamma_{U_t}(\Omega_{U_t|k}^{q+1}) \to \ldots) \ ,
$$
with $\nabla_t (\omega) = d \omega + dg_t \wedge \omega $ for a
local section $\omega $ of $\Omega_{U_t|k}^q $.

Now, if we have a global rapid decay chain $c \otimes e^{g_t} \in
\Gamma(\widetilde{X_t},\Crd_p(e^{-g_t})) $ with respect to the dual bundle
$\OU e^{-g_t} $ and a meromorphic $p $-form $\omega $, then the integral
$\int_c \omega e^{g_t} $ converges because the rapid decay of $e^{g_t} $
along $c $ annihilates the moderate growth of the meromorphic $\omega
$.  Let $c_\tau $ denote the topological chain one gets by cutting off a
small tubular neighborhood with radius $\tau $ around the boundary
$\partial \Delta^p $ from the given topological chain $c $. Then, for
$c \otimes e^{g_t} \in \Crd_p(e^{-g_t})(\widetilde{X_t}) $ and $\eta \in
\Omega^{p-1}(\sD) $ a meromorphic $(p-1) $-form, we have the {\bf 'limit
  Stokes formula'}
$$
\Int_c (\nabla_t \eta) e^{g_t} = \lim_{\tau \to 0} \Int_{c_\tau} (\nabla_t \eta)
e^{g_t} = \lim_{\tau \to 0} \Int_{\partial c_\tau} \eta \, e^{g_t} =
\Int_{\partial c-D_t} \eta \, e^{g_t} \,
$$
where in the last step we used that by the given growth/decay
conditions the integral over the faces of $\partial c_\tau $ 'converging'
against the faces of $\partial c $ contained in $D_t $ vanishes.

The limit Stokes formula easily shows in the standard way that
integrating a closed differential form over a given rd-cycle (i.e.
with vanishing boundary value) only depends on the de Rham class of
the differential form and the rd-homology class of the cycle. Thus, we
have:
\begin{prop} \label{prop:perpair} Integration induces a well defined
  bilinear pairing
  \begin{equation} \label{eq:perfV} H^p_{dR}(U_t, e^{g_t}) \times
    \Hrd_p(U_t, e^{-g_t}) \lto \C \ , \ ([\omega], [c \otimes
    e^{g_t}]) \mapsto \Int_c \omega \, e^{g_t} \ ,
  \end{equation}
  which we call the {\bf period pairing} of $\OU e^{g_t} $.
\end{prop}

\subsection{Duality}

We will now prove the following duality theorem for the connection
$\OU e^{g_t} $ in arbitrary dimension $\dim U_t $.

\begin{theo} \label{thm:dual} The period pairing \eqref{eq:perfV} is a
  perfect pairing of finite dimensional vector spaces.
\end{theo}

Note that one has the analogous result for any connection $(E, \nabla)
$ with good formal structure in dimension $\dim U_t \le 2 $
(\cite{hien}, or \cite{b-e} for $\dim U_t=1 $ resp.). The proof given in
the following is a generalization of the one in \cite{hien} for the
special case of a connection $\OU e^{g_t} $ to arbitrary dimension
$\dim(X_t) \in \N $. For the readers's convenience, we give the outline
of the proof and also comment on the necessary changes for higher
dimensions.

\begin{proof}
The proof will rely on a
local duality on $\widetilde{X_t} $. To prepare for it, we will describe
the period pairing as a local pairing of complexes of sheaves on
$\widetilde{X_t} $. We recall the corresponding notions of \cite{hien}.

Let $\Amod $ denote the sheaf of functions on $\widetilde{X_t} $ which
are holomorphic on $\Uan \subset \widetilde{X_t} $ and of moderate
growth along $\pi^{-1}(D_t) $ and let $\AD $ denote the subsheaf of such
functions that are infinitely flat along $\widetilde{D_t} $, i.e. all of
whose partial derivatives of arbitrary order vanish on $\widetilde{D_t}
$.

These sheaves are both flat over $\pi^{-1}(\O_{X_t}) $ and we call
$$
\DRmod(e^{g_t}) := \Amod \otimes_{\pi^{-1}(\O_{X_t})}
\pi^{-1}(\DR(\O_{X_t}[ \ast D_t] e^{g_t}))
$$
the {\em moderate de Rham complex} of $\O_{U_t} e^{g_t} $ and
$$
\DRD(e^{g_t}) := \AD \otimes_{\pi^{-1}(\O_{X_t})}
\pi^{-1}(\DR(\O_{X_t}[\ast D_t] e^{g_t}))
$$
the {\em asymptotically flat de Rham complex}. Note that the first one
computes the meromorphic de Rham cohomology of the meromorphic
connection $\O_{X_t}[\ast D_t] e^{g_t} $ (cp.
\cite{sabbah1}, Corollaire 1.1.8).

Now, the wedge product obviously defines a pairing of complexes
of sheaves on $\widetilde{X_t} $:
\begin{equation} \label{eq:locdu} \DRmod(e^{g_t}) \otimes_\C \DRD(e^{-g_t})
  \to \DRD(\O_{X_t}, d) \ ,
\end{equation}
where the right hand side denotes the asymptotically flat de Rham
complex of the trivial connection. We will compare this pairing to the
periods pairing from above.

Before that, let us introduce the sheaf $\Dbrd{-s} $ of rapid decay
distributions on $\widetilde{X_t} $, the local section of which are
distributions
$$
\vi \in \Db_{\widetilde{X_t}}^{-s}(V) := \Db_{\widetilde{X_t}}^{-s}(V) :=
\Hom_{\rm cont}( \Gamma_c(V, \Etop{s}), \C)
$$
on the space $\Etop{s} $ of $C^\infty $ differential forms on
$\widetilde{X_t} $ of degree $s $ with compact support in $\widetilde{X_t}
$ satisfying the following condition: we choose coordinates $x_1,
\ldots, x_n $ on $X_t $ such that locally on $V $ one has $D_t=\{ x_1
\cdots x_k =0 \} $. Then we require that for any compact $K \subset V $
and any element $N \in \N^k $ there are $m \in \N $ and $C_{K,N} >0 $
such that for any test form $\eta $ with compact support in $K $ the
estimate
\begin{equation} \label{eq:rddistrib} |\vi(\eta)| \le C_{K,N} \sum_i
  \sup_{|\alpha| \le m} \sup_K \{ |x|^N |\partial^\alpha f_i| \}
\end{equation}
holds, where $\alpha $ runs over all multi-indices of degree less than
or equal to $m $ and $\partial^\alpha $ denotes the $\alpha $-fold
partial derivative of the coefficient functions $f_i $ of $\eta $ in
the chosen coordinates.

Integration along chains induces a pairing of complexes of sheaves
\begin{equation} \label{eq:integral} \DRmd{s}(e^{g_t}) \otimes
  \Crdtp{-r}(e^{-g_t}) \to \Dbrd{s-r}
\end{equation}
to be defined as follows: for any local $s $-form $\omega $ of
$\DRmod(e^{g_t}) $, any rapid decay chain $c \otimes e^{g_t} \in
\Gamma(V, \Crdtp{-r}(e^{-g_t})) $ and any test form $\eta \in \Gamma_c(V,
\Etop{p}) $ with $p=r-s $, the decay and growth assumptions ensure
that the integral
\begin{equation} \label{eq:integralex} \int_{c} \eta \wedge ( e^{g_t}
  \cdot \omega )
\end{equation}
converges and moreover, that the distribution that maps $\eta $ to
\eqref{eq:integralex} satisfies the condition \eqref{eq:rddistrib}
from above.

Note that $\Dbrd{-\bdot} $ is a fine resolution of
$\widetilde{\jmath}_! \C_{\Uan}[2d] $ where $\widetilde{\jmath}:\Uan
\hookrightarrow \widetilde{X_t} $ denotes the inclusion and
$d:=\dim_{\C}(X_t) $. It follows that $\Hy^0(\widetilde{X_t},
\Dbrd{\bdot}) \cong \C $ in a standard way. The resulting pairing
$H^p_{dR}(U_t, e^{g_t}) \times \Hrd_p(U_t, e^{-g_t}) \to \C $ induced by
\eqref{eq:integral} on cohomology in degree zero coincides with the
period pairing (Proposition \ref{prop:perpair}).

The proof of Theorem \ref{thm:dual} now splits into the following
intermediate steps (cp. with \cite{hien}):

\begin{claim} \label{claim:deg0} The complexes $\DRmod(e^{g_t}) $ and
  $\DRD(e^{-g_t}) $ have cohomology in degree zero only.
\end{claim}
\claimpf In both cases the argument relies on an existence result of a
certain linear partial differential system with moderate or rapidly
decaying coefficients. To be more precise, consider the local
situation at some point $x_0 \in D_t=\{ x_1 \cdots x_k=0 \} $ and let
$\vt \in \pi^{-1}(x_0) \simeq (S^1)^k $ be a direction in
$\widetilde{D_t} $ over $x_0 $. Then the complex of stalks at $\vt $ which we
have to consider is given as
\begin{equation} \label{eq:PDE} \ldots \lto \Big( \As{X_t}{? D_t}
  \otimes_{\pi^{-1}\O_{X_t}} \pi^{-1} \Omega_{X_t}^p \Big)_\vt
  \stackrel{\nabla_t}{\lto} \Big( \As{X_t}{? D_t} \otimes_{\pi^{-1}\O_{X_t}}
  \pi^{-1} \Omega_{X_t}^{p+1} \Big)_\vt \lto \ldots \ ,
\end{equation}
where ? stands for either $< $ or mod. If we describe $\Omega_{X_t}^p $ in
terms of the local basis $dx_I $ for all $I=\{ 1 \le i_1 < \ldots <
i_p \le d\} $, the covariant derivation $\nabla_t $ in degree $p $ reads
as
$$
\sum_{\# I=p} w_I \, dx_I \mapsto \sum_{\# J=p+1} \Big( \sum_{j \in J}
sgn_J(j) (Q_j w_{J \smallsetminus \{ j \}}) \Big) dx_J
$$
where
$$
Q_j u:= \frac{\partial}{\partial x_j} u \pm \frac{\partial
  g_t}{\partial x_j} \cdot u \ ,
$$
with the negative sign in case of $e^{g_t} $ and the positive in case of
$e^{-g_t} $. Also, we let $sgn_J(j)=(-1)^\nu $ for $J=\{ j_1 < \ldots <
j_{p+1} \} $ and $j_\nu=j $.

Now, let $\omega $ be a germ of a section of $\As{X_t}{? D_t} \otimes
\pi^{-1} \Omega_{X_t}^{p+1} $ written in the chosen coordinates as
$$
\omega = \sum_{\# J= p+1} w_J \, dx_J \ ,
$$
such that $\nabla_t \omega = 0 $. We have to find a $p $-form $\eta $
with appropriate growth condition such that $\nabla_t \eta= \omega $.

To this end, let $r \in \N $ be an integer such that $w_J=0 $ for all
$J $ with $J \cap \{ 1, \ldots, r-1 \} \neq \emptyset $ (which is an
empty condition for $r=1 $). We prove that we can find a $p $-form
$\eta $ with coefficients in $\As{X_t}{? D_t} $ such that
\begin{equation} \label{eq:induction} \big( \omega - \nabla_t \eta \big)
  \in \sum_{J \cap \{1, \ldots, r\}= \emptyset} \Ass{X_t}{\vt}{? D_t} \,
  dx_J \ .
\end{equation}
Successive application of this argument will prove the claim.

By assumption
\begin{equation} \label{eq:omegaK} 0=\nabla_t \omega =: \sum_{\#K=p+2}
  \Big( \sum_{k \in K} sgn_K(k) (Q_k w_{K \smallsetminus \{ k\}})
  \Big) dx_K \ .
\end{equation}
Taking $k<r $ and $r \in J $ and examining the summand of
\eqref{eq:omegaK} corresponding to the set $K:= \{s \} \cup J $ we see
that
\begin{equation}\label{eq:intgby}
  Q_s w_J = 0 \mbox{\quad for all such } s=1, \ldots, r-1 \ .
\end{equation}

Now, consider the system $(\Sigma_J) $ of partial differential
equations for the unknown function $u_J $, where $J $ is a fixed
subset $J \subset \{ r, \ldots, d\} $ of cardinality $p+1 $ with $r
\in J $:
$$
(\Sigma_J): \left\{
  \begin{array}{ll}
    Q_s u_J= 0 & \mbox{for all } s=1, \ldots, r-1 \\[0.3cm]
    Q_r u_J= w_J \ ,
  \end{array}\right.
$$
together with the integability assumption \eqref{eq:intgby}. Systems
of this type had been studied by Majima (\cite{majima}) before. See
for example \cite{sabbah2}, Appendix A for a presentation from which
follows that this system always has a solution in $\Ass{X_t}{\vt}{<D_t} $
if the coefficients of the system belong to the same space of
functions. In the case of moderate growth coefficients/solutions, the
assertion follows from \cite{hien}, Theorem A.1 (which is formulated
for dimension 2 only, but generalizes one-to-one to the case of
arbitrary dimension).

In each case, we can always find a solution $u_J \in \Ass{X_t}{\vt}{? D_t}
$ for any such $J \subset \{r, \ldots, d\} $ of cardinality $p+1 $ and
$r \in J $ and if we let
$$
\eta:= \sum_{J \text{ as above}} u_J \, dx_J \ ,
$$
we easily see that \eqref{eq:induction} is satisfied.\\
\qd{Claim 1}

\begin{claim} \label{claim:perfder} The local pairing \eqref{eq:locdu}
  is perfect in the derived sense (i.e. the induced morphism
$$
\DRmod(e^{g_t}) \to \bR \Hom_{\widetilde{X_t}} (\DRD(e^{-g_t}),
\widetilde{\jmath}_! \C)
$$
and the analogous one with $\DRD $ and $\DRmod $ exchanged are
isomorphisms).
\end{claim}
\claimpf According to Claim \ref{claim:deg0} both complexes have
cohomology in degree $0 $ only, i.e. we are reduced to look at the
pairing of sheaves
$$
\Smod \otimes \SD \to \widetilde{\jmath}_! \C_{U_t}
$$
with $\Smod:= \H^0(\DRmod(e^{g_t})) $ and $\SD:=\H^0(\DRD(e^{-g_t})) $.
Again, we consider the local situation $D_t=\{ x_1 \cdots x_k=0 \} $ and
$g_t(x)=x_1^{-m_1} \cdots x_k^{-m_k} u(x) $, which we restrict to a
samll enough open polysector $V \subset \widetilde{X_t}
\stackrel{\pi}{\to} X_t $. We define the Stokes multi-directions of $g_t $
along $D_t $ inside $V $ to be
$$
\Sigma^{D_t}_{g_t}:= St^{-1} \big( ( \frac{\pi}{2}, \frac{3 \pi}{2}) \big) \ ,
$$
where $St: \pi^{-1}(D_t) \cap V \to \R/ 2\pi \Z $ denotes the map
$$
St(r_i, \vt_i) := - \sum_{i=1}^k m_i \vt_i + \arg(u \circ \pi(r_i,
\vt_i)) \ .
$$
Let $V_{g_t}:= (V \smallsetminus \pi^{-1}(D_t)) \cup
\Sigma_{g_t}^{D_t} $, i.e. $V_{g_t} \cap \pi^{-1}(D_t) $ are the
directions in which $e^{g_t(x)} $ has rapid decay for $x $ radially
approaching $D_t $. If $\widetilde{\jmath}_{g_t}:V_{g_t}
\hookrightarrow V $ denotes the inclusion, one obviously has (possibly
after shrinking $V $):
\begin{equation} \label{eq:Valpha} \Smod|_V =
  \widetilde{\jmath}_{-g_t}(e^{-g_t(x)} \cdot \C_{U_t})|_V \mbox{\quad and
    \quad} \SD|_V = \widetilde{\jmath}_{g_t} ( e^{g_t(x)} \cdot
  \C_{U_t})|_V \ .
\end{equation}
Consequently we have the following commutative diagramm:
\begin{equation} \label{eq:locV}
  \begin{CD}
    \bR \cHom ( \SD, \widetilde{\jmath}_!
    \C_{X_t \smallsetminus D_t} )|_V @>>> \Smod|_V  \\
    @V{\cong}VV @VV{\cong}V \\
    \bR \cHom \big( (\widetilde{\jmath}_{-g_t})_! \C_{V \smallsetminus
      D_t} \, , \, \widetilde{\jmath}_! \C_{V \smallsetminus D_t} \big)|_V
    @>>> (\widetilde{\jmath}_{g_t})_! \, \C_{V_{g_t}}
  \end{CD}
\end{equation}
By the factorization $\widetilde{\jmath} = \widetilde{\jmath}_{-g_t}
\circ \widetilde{\iota}_{-g_t} $ with $\widetilde{\iota}_{-g_t}: V
\smallsetminus \pi^{-1}(D_t) \hookrightarrow V_{-g_t} $, we see that
\begin{multline*}
  \bR \cHom \big( (\widetilde{\jmath}_{-g_t})_!  \C_{V_{-g_t}} \, , \,
  \widetilde{\jmath}_! \C_{X_t \smallsetminus D_t} \big) \cong
  (\widetilde{\jmath}_{-g_t})_\ast \bR \cHom \big( \C_{V_{-g_t}} ,
  (\widetilde{\iota}_{-g_t})_! \C_{V \smallsetminus \pi^{-1}(D_t)} \big)  \\
  \cong (\widetilde{\jmath}_{-g_t})_\ast \cHom \big( \C_{V_{-g_t}} ,
  (\widetilde{\iota}_{-g_t})_! \C_{V \smallsetminus \pi^{-1}(D_t)} \big) =
  (\widetilde{\jmath}_{-g_t})_! \, \C_{V_{g_t}} \ ,
\end{multline*}
since $(V \smallsetminus V_{g_t}) \cap \pi^{-1}(D_t) $ coincides with the
closure of $V_{-g_t} \cap \pi^{-1}(D_t) $ inside $\pi^{-1}(D_t) $. Hence,
the bottom line of \eqref{eq:locV} is an isomorphism and thus
$$
\Smod \cong \bR \cHom_{\widetilde{X_t}} \big( \SD, \widetilde{\jmath}_!
\C_{X_t \smallsetminus D_t} \big)
$$
locally on $\widetilde{X_t} $ over an arbitrary point of $D_t $.
Interchanging
$\SD $ and $\Smod $ gives the analogous isomorphism.\\
\qd{Claim 2}

\begin{claim} \label{claim:crdtDRD} There is a canonical isomorphism
$$
\Crdt(e^{-g_t}) \cong \DRD(e^{-g_t})[2d]
$$
in the derived category $D^b(\C_{\widetilde{X_t}}) $, where
$d=\dim_\C(X_t) $.
\end{claim}
\claimpf Here, the arguments of \cite{hien}, Theorem 3.6, apply
directly, we therefore only briefly give the main line of the proof:
By Claim \ref{claim:deg0} we have
$$
\DRD(e^{-g_t}) \stackrel{\sim}{\longleftarrow} \H^0(\DRD(e^{-g_t}))=: \SD
\ .
$$
Additionally, there is a canonical morphism
$$
\Ctd^{-\bdot} \otimes \SD \lto \Crdt(e^{-g_t}) \ ,
$$
since for any open $V \subset \widetilde{X_t} $ and an asymptotically
flat section $\sigma \in \Gamma_V(\SD) $, the section $\sigma \in
\SD(V) \subset \widetilde{\jmath}_\ast (e^{g_t} \C_{U_t}) $ is rapidly
decaying along any chain in $c \in \Ctd^{-\bdot}(V) $, hence $c
\otimes \sigma \in \Gamma(V,\Crdt(e^{-g_t})) $. We prove that this
morphism is a quasi-isomorphism. Restricted to $U_t \subset
\widetilde{X_t} $, this is obvious.

Next, consider a small open polysector $V $ around some $\vt \in
\pi^{-1}(x_0) $ with $x_0 \in D_t $. Let $\Sigma^{D_t}_{\pm g_t} \subset
\pi^{-1}(x_0) $ denote the set of Stokes-directions of $e^{\pm g_t} $ as
in the proof of Claim \ref{claim:deg0}.

If $\vt \in \Sigma^{D_t}_{-g_t} $ then for $V $ being a small enough
polysector, we have $\SD |_V = \widetilde{\jmath}_! \big( e^{g_t} \C_{U_t}
\big) |_V $. For a smooth topological chain $c $ in $\widetilde{X_t} $,
the local section $e^{g_t} $ will not have rapid decay along $c $ in $V
$ as required by the definition unless the chain does not meet
$\widetilde{D_t} \cap V $. Hence
$$
\Crdt(e^{-g_t})|_V = \Ctd^{-\bdot} \otimes \SD|_V \ .
$$

If $\vt \in \Sigma^{D_t}_{g_t} $, we can assume that $V $ is an open
polysector such that all the arguments of points in $V $ are contained
in $\Sigma^{D_t}_{g_t} $. Then $\SD|_V \cong \widetilde{\jmath}_\ast
(e^{g_t} \C_V) $.  Similarly, all twisted chains $c \otimes e^{g_t} $
will have rapid decay inside $V $ and again both complexes considered
are equal to $\Ctd^{-\bdot} \otimes \widetilde{\jmath}_\ast (e^{g_t}
\C_{U_t}) $.

Finally, if $\vt $ separates the Stokes regions of $e^{g_t} $ and $e^{-g_t}
$, we have with the notation of \eqref{eq:Valpha}:
$$
\SD|_V \cong (\widetilde{\jmath}_{g_t})_! (e^{g_t} \C_{U_t})|_V \ .
$$
The subspace $V_{g_t} $ is characterized by the property that $V_{g_t}
\cap \widetilde{D_t} $ consists of those directions along which
$e^{g_t(x)} $ has rapid decay for $x $ approaching $\widetilde{D_t} $. In
particular, $c \otimes e^{g_t} $ is a rapid decay chain on $V $ if and
only if the topological chain $c $ in $\widetilde{X_t} $ approaches
$\widetilde{D_t} \cap V $ in $V_{g_t} $ at most. Hence
$$
\Crdt(e^{-g_t})|_V = \Ctd^{-\bdot} \otimes ( \widetilde{\jmath}_{g_t} )_!
(e^{g_t} \C_{U_t}) = \Ctd^{-\bdot} \otimes \SD|_V \ .
$$

In summary, we have the following composition of quasi-isomorphisms
$$
\SD[2d] \stackrel{\simeq}{\lto} \C_{\widetilde{X_t}}[2d] \otimes \SD
\stackrel{\simeq}{\lto} \Ctd^{-\bdot} \otimes \SD
\stackrel{\simeq}{\lto} \Crdt(e^{-g_t}) \ .
$$
\qd{Claim 3}

\noindent
In order to prove Theorem \ref{thm:dual}, consider the resolution
$$
\AD \otimes_{\pi^{-1}\O_{X_t}} \pi^{-1} \Omega_{X_t}^r \hookrightarrow
\big( \mathcal{P}_{\widetilde{X_t}}^{<D_t} \otimes_{\pi^{-1} C_{X_t}^\infty}
\pi^{-1} \Epq{r}{\bdot}, \dolb \big) \ ,
$$
where $\mathcal{P}_{\widetilde{X_t}}^{<D_t} $ denotes the sheaf of
$C^\infty $-functions flat at $\pi^{-1}(D_t) $ and $\Epq{r}{s} $ denotes
the sheaf of $C^\infty $ forms on $X_t $ of degree $(r,s) $. We then
have the bicomplex
$$
\Rd^{\bdot, \bdot}:= \big( \mathcal{P}_{\widetilde{X_t}}^{<D_t}
\otimes_{\pi^{-1}C_{X_t}^\infty} \pi^{-1} \Epq{\bdot}{\bdot} \, ,
\partial, \dolb \big) \ ,
$$
whose total complex $\Rd^{\bdot} $ computes is a fine resolution of
$\widetilde{\jmath}_! \C_{U_t} $.

According to Claim \ref{claim:deg0}, the local duality pairing
\eqref{eq:locdu} reads as $\SD \otimes \Smod \to \widetilde{\jmath}_!
\C_{\Uan} $. We now have a canonical quasi-isomorphism
$$
\beta: \Ctd^{-\bdot} \otimes \Rd^{\bdot} \stackrel{\simeq}{\lto}
\Dbrd{-\bdot}
$$
of complexes mapping an element $c \otimes \rho \in \Ctd^{-r} \otimes
\Rd^{s}(V) $ of the left hand side over some open $V \subset
\widetilde{X_t} $ to the distribution given by $\eta \mapsto \int_c \eta
\wedge \rho $ for a test form $\eta $ with compact support in $V $.
The latter represents a local section of $\Dbrd{s-r} $ since $\rho $
is rapidly decaying. Additionally, we have to consider the morphism
$$
\gamma:\Crdt(e^{-g_t}) \otimes \Smod \to \Dbrd{-\bdot} \ , \ (c \otimes
e^{g_t}) \otimes \sigma \mapsto ( \eta \mapsto \int_c e^{g_t} \cdot
\sigma \cdot \eta ) \ ,
$$
which induces the period pairing after taking cohomology.

These morphisms fit into the following commutative diagramm
\begin{equation}\label{eq:finalcommdiag}
  \begin{CD}
    \SD[2d] \otimes \Smod @>>> \widetilde{\jmath}_! \C_{\Uan}[2d] \\
    @V{\simeq}VV @VV{\simeq}V \\
    \Ctd^{-\bdot} \otimes \SD \otimes \Smod @>>> \Ctd^{-\bdot}
    \otimes \Rd^{\bdot} \\
    @V{\simeq}V{\mbox{\scriptsize Claim \ref{claim:crdtDRD}}}V
    @V{\simeq}V{\beta}V \\
    \Crdt(e^{-g_t}) \otimes \Smod @>{\gamma}>> \Dbrd{-\bdot} \ .
  \end{CD}
\end{equation}

The top row is a perfect duality in the derived sense by Claim
\ref{claim:perfder}. Applying Poincar\'e-Verdier duality thus induces
an isomorphism
\begin{multline}\label{eq:PV}
  \bR \Gamma_{\widetilde{X_t}} \DRD(e^{-g_t})[2d] \cong \bR
  \Gamma_{\widetilde{X_t}} \bR \cHom_{\widetilde{X_t}}( \DRmod(e^{g_t}),
  \widetilde{\jmath}_!\C_{U_t}[2d])
  \cong \\
  \cong \Hom_\C^{\bdot} (\bR \Gamma_{\widetilde{X_t}} \DRmod(e^{g_t}), \C )
  \ .
\end{multline}
The commutativity of \eqref{eq:finalcommdiag} shows that the morphism
induced by the period pairing i.e. by $\gamma $, namely
$$
\bR \Gamma_{\widetilde{X_t}} \Crdt(e^{-g_t}) \lto \Hom_\C^{\bdot} (\bR
\Gamma_{\widetilde{X_t}} \DRmod(e^{g_t}), \C ) \ ,
$$
is an isomorphism. The same holds after exchanging $\Crdt $ and
$\DRmod $. Taking $p $-th cohomology therefore yields the perfect
period pairing
$$
\Hdr^p(U_t, e^{g_t}) \otimes_\C \Hrd_p(U_t, e^{-g_t}) \to \C
$$
completing the proof of Theorem \ref{thm:dual}.\\
\end{proof}

We will now study the exponential Gau\ss-Manin systems and apply these
considerations to describe the holomorphic solutions of the latter.

\section{Holomorphic solutions of exponential Gau\ss-Manin
  systems}

Consider the situation described in the introduction. Let $f,g:U\to\mathbb{A}^1$ be two algebraically independant
regular functions on the smooth affine variety $U$.
We denote by $X$ a smooth projective compactification of $U$, on which
$f,g:U\to\mathbb{A}^1$ extend to functions $F,G:X\to\mathbb{P}^1$. Let
$D:=X\setminus U$ and assume that it is a normal crossing divisor.

\subsection{Definition}
We define an exponential Gau\ss-Manin systems as being the
cohomology sheaves of the direct image by $f$ of the module of
exponential type $\mathcal{O}_Ue^g$:
$$\mathcal{G}^k:=\mathcal{H}^kf_+(\mathcal{O}_Ue^g).$$

Let $\Omega_U^k$ be the sheaf of algebraic differential forms of degree
$k$ on $U$.  Using the canonical factorisation of $f$ into a closed
embedding and a projection, we prove as in \cite{disa} that:

\begin{prop}\label{definition}
  $f_+(\mathcal{O}_Ue^g)$ is isomorphic to 
  $f_{\ast}(\Omega_U^{\bullet+n}\otimes_{\mathbb{C}}\mathbb{C}[\partial_t])$,
  where the differential $\nabla$ on
  $\Omega_U^{\bullet+n}\otimes_{\mathbb{C}}\mathbb{C}[\partial_t]$ is
  given by
$$\nabla(w\otimes\partial_t^i)=dw\otimes\partial_t^i+dg\wedge
w\otimes\partial_t^i-df\wedge w\otimes\partial_t^{i+1}.
$$
Its cohomology modules $\mathcal{G}^k$ are equipped with a structure
of $\mathcal{D}_{\mathbb{A}^1}$-module given by:
$$\begin{array}{l}
  t[w\partial_t^i]=[fw\partial_t^i-iw\partial_t^{i-1}],\\
  \partial_t[w\partial_t^i]=[w\partial_t^{i+1}].
\end{array}$$
They are holonomic $\mathcal{D}_{\mathbb{A}^1}$-modules with irregular
singularities on the affine line $\mathbb{A}^1$ and at infinity.
\end{prop}

\begin{rema}
The interpretation of the regular singularities of $\mathcal{G}^k$ in terms of invariants associated to $f$ and $g$ is still unknown. Nevertheless, in \cite{roucairol} and \cite{roucairol2}, a complete description of the irregular singularities is given. Let $\Delta_k\subset\mathbb{A}^1\times\mathbb{A}^1$ be the singular support of the holonomic $\mathcal{D}_{\mathbb{A}^1\times\mathbb{A}^1}$-module 
  $\mathcal{H}^k(f,g)_+(\mathcal{O}_U)$. Its closure $\overline{\Delta_k}$ in $\mathbb{P}^1\times\mathbb{P}^1$ is included in the projective variety $\overline{\Delta}$, such that $(F,G):(F,G)^{-1}((\mathbb{P}^1\times\mathbb{P}^1)\setminus\overline{\Delta})\to(\mathbb{P}^1\times\mathbb{P}^1)\setminus\overline{\Delta}$ is a locally trivial fibration. $t_0\in\mathbb{A}^1\cup\{\infty\}$ is an
  irregular singularity of $\mathcal{G}^k$ if and only if
  the germ of
  $\overline{\Delta_k}$ at $(t_0,\infty)$ is not empty and not included in $\{t_0\}\times\mathbb{P}^1$.
\end{rema}

\subsection{Results on the exponential Gau\ss-Manin systems}
\subsubsection{Exponential Gau\ss-Manin connection and its sheaf of flat sections}

Let $\Sigma_1$ be a finite subset of $\mathbb{A}^1$ such that:
\begin{itemize}
\item $f:f^{-1}(\mathbb{A}^1\setminus\Sigma_1)\to\mathbb{A}^1\setminus\Sigma_1$ is a locally trivial fibration. For any
$t\in\mathbb{A}^1\setminus\Sigma_1$, $f^{-1}(t)$ is a smooth affine
variety. We will write $g_t$ instead of $g_{|f^{-1}(t)}$.
\item $\mathcal{G}^k_{|\mathbb{A}^1\setminus\Sigma_1}$ is locally free, hence a flat
connection, for any $k\in\mathbb{Z}$ ($\mathcal{G}^k$ is holonomic).
\end{itemize}

Let $(\Gamma_{f^{-1}(\mathbb{A}^1\setminus\Sigma_1)}(\Omega^{\bullet+n-1}_{U|\mathbb{A}^1}),\nabla)$ be the complex of relative algebraic differential forms on $f^{-1}(\mathbb{A}^1\setminus\Sigma_1)$ with respect to $f$ equipped with the differential defined by $\nabla(w)=dw+dg\wedge w$. 

\begin{prop}\label{horizontal}
\begin{enumerate}
\item $\Gamma(\mathbb{A}^1\setminus\Sigma_1,\mathcal{G}^k)$ is a $\mathcal{D}(\mathbb{A}^1\setminus\Sigma_1)$-module isomorphic to $H^k(\Gamma_{f^{-1}(\mathbb{A}^1\setminus\Sigma_1)}(\Omega^{\bullet+n-1}_{U|\mathbb{A}^1}),\nabla)$. The action of $\mathcal{D}(\mathbb{A}^1\setminus\Sigma_1)$ on a form $[w]\in H^k(\Gamma_{f^{-1}(\mathbb{A}^1\setminus\Sigma_1)}(\Omega^{\bullet+n-1}_{U|\mathbb{A}^1}),\nabla)$ is given by:
$$\begin{array}{lll}
  t.[w]&=&[fw],\\
  \partial_t.[w]&=&[\eta],
\end{array}$$
where $[\eta]\in H^k(\Gamma_{f^{-1}(\mathbb{A}^1\setminus\Sigma_1)}(\Omega^{\bullet+n-1}_{U|\mathbb{A}^1}),\nabla)$ satisfies $dw+dg\wedge w=df\wedge\eta$ (which comes from the condition $\nabla(w)=0$).
\item The sheaf of flat sections of $\mathcal{G}^k$ is a local system on $\mathbb{A}^1\setminus\Sigma_1$ and the stalk at the point $t$ is $\Hdr^{k+n-1}(f^{-1}(t),e^{g_t})$. 

This local system on $\mathbb{A}^1\setminus\Sigma_1$ will be denoted by $\widetilde{\Hdr^{k+n-1}}$ and we will consider $(\Hdr^{k+n-1},\nabla)$ its associated flat vector bundle.
\end{enumerate}
\end{prop}

\noindent
The action of $t$, $\partial_t$ does not depend on the class of $w$ in
$H^k(\Gamma_{f^{-1}(\mathbb{A}^1\setminus\Sigma_1)}(\Omega^{\bullet+n-1}_{U|\mathbb{A}^1}),\nabla)$.
Indeed, if $[w]=0$, we have $w=df\wedge\alpha+d\beta+dg\wedge \beta$, where $\alpha$ and $\beta$ are some forms. Then,
\begin{itemize}
\item $t.[w]=[fw]=[df\wedge(f\alpha-\beta)+d(f\beta)+dg\wedge f\beta]=0$.
\item As $dw+dg\wedge w=-df\wedge(d\alpha+dg\wedge\alpha)$, we have $\partial_t.[w]=-[d\alpha+dg\wedge\alpha]=0$.
\end{itemize}

\begin{proof}

\noindent
\begin{enumerate}
\item As $f_{|f^{-1}(\mathbb{A}^1\setminus\Sigma_1)}$ is smooth,
$f_+(\mathcal{O}_Ue^g)_{|\mathbb{A}^1\setminus\Sigma_1}=\bR f_{\ast}(\Omega^{\bullet+n-1}_{U|\mathbb{A}^1})_{|\mathbb{A}^1\setminus\Sigma_1}$
(cf. Pro\-po\-sition 1.4 in \cite{dmss}). Then, $\Gamma(\mathbb{A}^1\setminus\Sigma_1,\mathcal{G}^k)=\Gamma(\mathbb{A}^1\setminus\Sigma_1,\mathcal{H}^k\bR f_{\ast}(\Omega^{\bullet+n-1}_{U|\mathbb{A}^1}))$. But we have the spectral sequence:
$$H^i(\mathbb{A}^1\setminus\Sigma_1,\mathcal{H}^{k-i}\bR f_{\ast}(\Omega^{\bullet+n-1}_{U|\mathbb{A}^1}))\Longrightarrow\mathbb{H}^k(\mathbb{A}^1\setminus\Sigma_1,\bR f_{\ast}(\Omega^{\bullet+n-1}_{U|\mathbb{A}^1})).$$
As $\mathbb{A}^1\setminus\Sigma_1$ is affine and $\mathcal{H}^{k-i}\bR f_{\ast}(\Omega^{\bullet+n-1}_{U|\mathbb{A}^1})$ is coherent, this spectral sequence degenerates and we have:
$$\Gamma(\mathbb{A}^1\setminus\Sigma_1,\mathcal{H}^k\bR f_{\ast}(\Omega^{\bullet+n-1}_{U|\mathbb{A}^1}))=\mathbb{H}^k(\mathbb{A}^1\setminus\Sigma_1,\bR f_{\ast}(\Omega^{\bullet+n-1}_{U|\mathbb{A}^1})).$$

Then,
$$\begin{array}{lll}
\Gamma(\mathbb{A}^1\setminus\Sigma_1,\mathcal{G}^k)
&=&\mathbb{H}^k(\mathbb{A}^1\setminus\Sigma_1,\bR f_{\ast}(\Omega^{\bullet+n-1}_{U|\mathbb{A}^1})),\\
&=&\mathbb{H}^k(f^{-1}(\mathbb{A}^1\setminus\Sigma_1),\Omega^{\bullet+n-1}_{U|\mathbb{A}^1}),\\
&=&H^k(\Gamma_{f^{-1}(\mathbb{A}^1\setminus\Sigma_1)}(\Omega^{\bullet+n-1}_{U|\mathbb{A}^1}),\nabla).
\end{array}$$

The last isomorphism comes from the coherence of $\Omega^{\bullet+n-1}_{U|\mathbb{A}^1~|f^{-1}(\mathbb{A}^1\setminus\Sigma_1)}$ and the fact that $f^{-1}(\mathbb{A}^1\setminus\Sigma_1)$ is affine.

Concerning the action of $t$ and $\partial_t$ on $H^k(\Gamma_{f^{-1}(\mathbb{A}^1\setminus\Sigma_1)}(\Omega^{\bullet+n-1}_{U|\mathbb{A}^1}),\nabla)$, we have to explicite the previous isomorphism. In fact, we can prove that it is given by:

$$\begin{array}{cccc}
\Psi:&H^k(\Gamma_{f^{-1}(\mathbb{A}^1\setminus\Sigma_1)}(\Omega^{\bullet+n-1}_{U|\mathbb{A}^1}),\nabla)&\to&\Gamma(\mathbb{A}^1\setminus\Sigma_1,\mathcal{G}^k))\\
&[w]&\mapsto&[-df\wedge w\otimes 1].
\end{array}$$ 

The injectivity and surjectivity of $\Psi$ is proved with similar calculations as in \cite{pham} (p.160) and is left to the reader.
It is a morphism of $\mathcal{D}(\mathbb{A}^1\setminus\Sigma_1)$-modules because:
\begin{itemize}
\item[$\bullet$] $\psi(t[w])=\psi([fw])=[-fdf\wedge w\otimes 1]=t.[-df\wedge w\otimes 1]=t.\psi([w])$.
\item[$\bullet$] $\psi(\partial_t[w])=\psi([\eta])=[-df\wedge\eta\otimes 1]=[-(dw+dg\wedge w)\otimes 1]=[-df\wedge w\otimes\partial_t]=\partial_t.[-df\wedge w\otimes 1]=\partial_t.\psi([w])$.
\end{itemize}

\item We recall the proof given in \cite{maaref}. As $\mathcal{G}^k$ is a flat connection on
$\mathbb{A}^1\setminus\Sigma_1$, the sheaf of flat sections of $\mathcal{G}^k$ is a local system on $\mathbb{A}^1\setminus\Sigma_1$ and the stalk at the point $t$ is 
$i^+_t\mathcal{G}^k=\mathcal{H}^ki_t^+f_+(\mathcal{O}_Ue^g)$. Now,
according to the base change formula (cf. Theorem 8.4 p. 267 in
\cite{borel}) and as $f_{|f^{-1}(\mathbb{A}^1\setminus\Sigma_1)}$ is
smooth,
$$
\begin{array}{lll}
  i_t^+\mathcal{G}^k&=&\mathcal{H}^kf_{t+}(\mathcal{O}_{f^{-1}(t)}e^{g_t}), 
\text{ with }f_t:=f_{|f^{-1}(t)}:f^{-1}(t) \to \{t\}, \\
  &=&\Hdr^{k+n-1}(f^{-1}(t),e^{g_t}). 
\end{array}$$
\end{enumerate}
\end{proof}

\subsubsection{Analytification of $\mathcal{G}^k$}
We will differentiate the analytic and the algebraic settings by writing an exponent $~^{\an}$ in the analytic situation.

Let
$(\mathcal{G}^k)^{\an}:=
\mathcal{O}_{\mathbb{C}}^{\an}\otimes_{\mathcal{O}_{\mathbb{A}^1}}\mathcal{G}^k$
be the analytification of $\mathcal{G}^k$.

From Proposition \ref{horizontal} we easily deduce:

\begin{prop}\label{prop:analytG}
  Let $V$ be an open subset of $\mathbb{C}\setminus\Sigma_1$.
  \begin{enumerate}
  \item $(\mathcal{G}^k)^{\an}(V)$ is isomorphic to
    the space of holomorphic sections of $(\Hdr^{k+n-1},\nabla)$ on $V$.
  \item $(\mathcal{G}^k)^{\an}(V)$ is isomorphic to $\mathcal{O}^{\an}(V)
    \otimes_{\mathbb{C}[t]}
    H^k(\Gamma_{f^{-1}(\mathbb{A}^1\setminus\Sigma_1)}(\Omega^{\bullet+n-1}_{U|\mathbb{A}^1}),\nabla)$; the action of $\mathcal{D}^{\an}(V)$ on $a\otimes[w]\in \mathcal{O}^{\an}(V)\otimes_{\mathbb{C}[t]}H^k(\Gamma_{f^{-1}(\mathbb{A}^1\setminus\Sigma_1)}(\Omega^{\bullet+n-1}_{U|\mathbb{A}^1}),\nabla)$ is given by:
        $$\begin{array}{lll}
  h.(a\otimes[w])&=&(ha)\otimes[w],~\text{where }h\in\mathcal{O}^{\an}(V),\\
  \partial_t.(a\otimes[w])&=&\frac{\partial a}{\partial t}\otimes[w]+a\otimes[\eta],
\end{array}$$
where $[\eta]\in H^k(\Gamma_{f^{-1}(\mathbb{A}^1\setminus\Sigma_1)}(\Omega^{\bullet+n-1}_{U|\mathbb{A}^1}),\nabla)$ satisfies $dw+dg\wedge w=df\wedge\eta$.

  \end{enumerate}
\end{prop}

\subsection{Family of rapid decay cycles}

\begin{prop}\label{hrd flat}
  There exists a finite subset $\Sigma_2$ of $\mathbb{C}$ such that
  $$\underset{t\in\mathbb{C}\setminus\Sigma_2}{\bigcup}\Hrd_k(f^{-1}(t),e^{-g_t})\to\mathbb{C}\setminus\Sigma_2$$
  is a flat vector bundle. We will denote $\widetilde{\Hrd_k}$ its associated local system, it being the sheaf of family of cycles in $\Hrd_k(f^{-1}(t),e^{-g_t})$ which depend continuously on $t$.
\end{prop}

Let $\mathcal{W}$ be a Whitney stratification of $X$ such that
$G^{-1}(\infty)$, $D$ and $U$ are union of strata (cf. Theorem (2.2)
in \cite{verdier2}). The central argument of this proof is the
following lemma.

\begin{lemm}\label{fibration}
  There exists a finite subset $\Sigma_2$ of $\mathbb{C}$ such that
  $F:F^{-1}(\mathbb{C}\setminus\Sigma_2)\to\mathbb{C}\setminus\Sigma_2$
  is a locally trivial fibration which respect the stratification on
  $F^{-1}(\mathbb{C}\setminus\Sigma_2)$ induced by $\mathcal{W}$.

\renewcommand{\labelenumi}{{\rm \roman{enumi})}}
  \begin{enumerate}
  \item This means that for any $t_0\in\mathbb{C}\setminus\Sigma_2$,
    there exists a small disc $S\subset\mathbb{C}\setminus\Sigma_2$
    centered at $t_0$ and a homeomorphism $\phi:F^{-1}(t_0)\times S\to
    F^{-1}(S)$ which respect the stratification induced by
    $\mathcal{W}$ and such that
$$\forall (x,t)\in F^{-1}(t_0)\times S,~F\circ\phi(x,t)=t.$$

\item Moreover, we can assume that there exists a real positive number
  $R$ big enough such that
$$\forall (x,t)\in \big(f^{-1}(t_0)\cap g^{-1}(\{|\rho|>R\})\big)\times S,~g\circ\phi(x,t)=g(x).$$
\end{enumerate}
\end{lemm}

\noindent
\prf{Lemma \ref{fibration}} 
The first point of this Lemma is the first
isotopy Lemma. Here, we have to adapt its proof in the way of checking
that there exists a local trivialisation such that the second point is
true. If we follow the proof of the first isotopy Lemma given in
\cite{verdier2} (cf. Theorem (4.14)), this can be done by constructing
some well-chosen rugous tangent vector fields on
$(F^{-1}(\mathbb{C}\setminus\Sigma_2),\mathcal{W})$. We just give a
sketch of the proof.

According to the Sard's Theorem (cf. Theorem (3.3) in
\cite{verdier2}), and the fact that $f$ and $g$ are algebraically
independant, we know that
\begin{enumerate}
\item there exists a finite subset $\Sigma_0$ of $\mathbb{C}$ such that
  $F:F^{-1}(\mathbb{C}\setminus\Sigma_0)\to\mathbb{C}\setminus\Sigma_0$
  is a morphism which is transverse to the stratification on
  $F^{-1}(\mathbb{C}\setminus\Sigma_0)$ induced by $\mathcal{W}$;
\item there exists a projective variety $\overline{\Delta}$ in
  $\mathbb{P}^1\times\mathbb{P}^1$ of dimension at most one such that
  $(F,G):(F,G)^{-1}(\mathbb{P}^1\times\mathbb{P}^1\setminus\overline{\Delta})\to\mathbb{P}^1\times\mathbb{P}^1\setminus\overline{\Delta}$
  is a morphism which is transverse to the stratification on
  $(F,G)^{-1}(\mathbb{P}^1\times\mathbb{P}^1\setminus\overline{\Delta})$
  induced by $\mathcal{W}$.
\end{enumerate}

We define $\Sigma_2$ as the union of $\Sigma_0$ and the set of
$t_0\in\mathbb{C}$ such that $(t_0,\infty)$ belongs to the closure of
$\overline{\Delta}\setminus(\mathbb{P}^1\times\{\infty\})$ in
$\mathbb{P}^1\times\mathbb{P}^1$.

Let $t_0\in\mathbb{C}\setminus\Sigma_2$. Then,
$\mathbb{C}\setminus\Sigma_2$ is isomorphic to $I_1\times I_2$ in a
neighbourhood of $t_0=(t_0^1,t_0^2)$, where $I_1$ and $I_2$ are open
intervals in $\mathbb{R}$. Let $R>0$ such that $\big((I_1\times
I_2)\times\{|\rho|>R\}\big)\cap\overline{\Delta}=\emptyset$.

If $t_i$ is a coordinate on $I_i$, $i=1,2$, we consider $\eta^i$
(resp. $\widetilde{\eta}^i$) the vector field
$\frac{\partial}{\partial t_i}$ on $I_1\times I_2$ (resp. $(I_1\times
I_2)\times\{|\rho|>R\}\subset\overline{\Delta}$). According to
Proposition (4.6) in \cite{verdier2}, there exists a rugous tangent
vector field $\xi^i$ on $(F^{-1}(I_1\times I_2),\mathcal{W})$ such
that

\begin{enumerate}
\item for all $x\in F^{-1}(I_1\times I_2)$,
  $TF(\xi^i_x)=\eta^i_{F(x)}$;
\item for all $x\in (f,g)^{-1}((I_1\times I_2)\times\{|\rho|>R\})$,
  $T(f,g)(\xi^i_x)=\widetilde{\eta}^i_{(f,g)(x)}$.
\end{enumerate}

Let $\phi^i:U^i\subset F^{-1}(I_1\times I_2)\times\mathbb{R}\to
F^{-1}(I_1\times I_2)$ be the flow of $\xi^i$, where $U^i$ is an open
neighbourhood in $F^{-1}(I_1\times I_2)\times\mathbb{R}$ of
$F^{-1}(I_1\times I_2)\times\{t_0^i\}$ (cf. Proposition (4.8) in
\cite{verdier2}).

As $F$ is a proper map, there exist $\widetilde{I_1}\subset I_1$ and
$\widetilde{I_2}\subset I_2$ two open intervals such that
$$
\begin{array}{lll}
  F^{-1}(t_0)\times(\widetilde{I_1}\times\widetilde{I_2}) & \to &
  F^{-1}(\widetilde{I_1}\times\widetilde{I_2})\\
  (x,(t_1,t_2))&\mapsto&\phi^2(\phi^1(x,t_1),t_2)
\end{array}$$ 
is a homeomorphism which satisfies properties i) and ii) of Lemma
\ref{fibration}.\\
\qed

\vspace{0.5cm}

\noindent
\prf{Proposition \ref{hrd flat}} We recall that a cycle in
$\Hrd_k(f^{-1}(t),e^{-g_t})$ is a relative cycle in
$H_k(F^{-1}(t),F^{-1}(t)\cap D,\mathbb{C})$ with a rapid decay
condition on $e^{g_t}$.  Lemma \ref{fibration} i) implies that
$$\pi:\underset{t\in\mathbb{C}\setminus\Sigma_2}{\bigcup}H_k(F^{-1}(t),F^{-1}(t)\cap D,\mathbb{C})\to\mathbb{C}\setminus\Sigma_2$$ 
is a flat vector bundle.
If $t_0\in\mathbb{C}\setminus\Sigma_2$ and $(S,\phi)$ is as in Lemma \ref{fibration}, a trivialisation
$$T:\underset{t\in S}{\bigcup}H_k(F^{-1}(t),F^{-1}(t)\cap D,\mathbb{C})\to H_k(F^{-1}(t_0),F^{-1}(t_0)\cap D,\mathbb{C})\times S$$
is given by $T(c_t)=(p_1\circ\phi^{-1}\circ c_t,t)$ for any $c_t\in H_k(F^{-1}(t),F^{-1}(t)\cap D,\mathbb{C})$, where
$p_1:F^{-1}(t_0)\times S\to F^{-1}(t_0)$ is the canonical projection. Moreover, for any $t\in S$, the restriction  $T_{tt_0}:H_k(F^{-1}(t),F^{-1}(t)\cap D,\mathbb{C})\to H_k(F^{-1}(t_0),F^{-1}(t_0)\cap D,\mathbb{C})$ of $T$ is a linear isomorphism.
Indeed, the inclusion of the pair $(F^{-1}(t),F^{-1}(t)\cap D)$ in $(F^{-1}(S),F^{-1}(S)\cap D)$ induces an isomorphism $i_t:H_k(F^{-1}(t),F^{-1}(t)\cap D,\mathbb{C})\to H_k(F^{-1}(S),F^{-1}(S)\cap D,\mathbb{C})$ and $i_t^{-1}$ is homotope to $p_1\circ\phi^{-1}$. Then, we have the following diagram:

$$\xymatrix{
&H_k(F^{-1}(t_0),F^{-1}(t_0)\cap D,\mathbb{C})\\
H_k(F^{-1}(S),F^{-1}(S)\cap D,\mathbb{C})\ar[ur]^{p_1\circ\phi^{-1}}_{\simeq}&\\
&H_k(F^{-1}(t),F^{-1}(t)\cap D,\mathbb{C})\ar[uu]_{T_{tt_0}}^{\simeq}\ar[ul]^{i_t}_{\simeq},
}
$$
which proves that $T_{tt_0}$ is an isomorphism. 
This gives the structure of vector bundle. Now, to show that it is flat, we just have to remark that the changes of trivialisations correspond to changes of basis in $H_k(F^{-1}(S),F^{-1}(S)\cap D,\mathbb{C})$, for some open set $S$ and thus, are constant.

Then, we have to verify that if we require the
rapid decay condition, we obtain a subbundle of this bundle. It is sufficient to prove that the rapid decay condition is
stable with respect to the isomorphism $T_{tt_0}$. In the following,
we will denote by $\psi_{tt_0}:F^{-1}(t)\to F^{-1}(t_0)$, the
homeomorphism $p_1\circ\phi^{-1}$.

Let $c_t\otimes e^g$ be a rapid decay chain in $F^{-1}(t)$. We want
to prove that $T_{tt_0}(c_t)\otimes e^g$ is a rapid decay chain in
$F^{-1}(t_0)$. Let $y\in T_{tt_0}(c_t)(\Delta^p)\cap D$. Then there
exists $\widetilde{y}\in c_t(\Delta^p)\cap D$ such that
$\psi_{tt_0}(\widetilde{y})=y$. As $e^g$ has rapid decay on $c_t$,
$G(\widetilde{y})=\infty$ and as $\psi_{tt_0}$ respects
$G^{-1}(\infty)$, we have also $G(y)=\infty$.

But according to Lemma \ref{fibration} ii), there exists a
neighbourhood $V$ of $y$ such that for all $z\in f^{-1}(t_0)\cap V$,
$g\circ\psi^{-1}_{tt_0}(z)=g(z)$. Then, if $e^g$ has rapid decay on
$c_t$ it has rapid decay on $T_{tt_0}(c_t)$. The converse can be
proved similary.\\
\qed

\subsection{Holomorphic solutions of the exponential Gauss-Manin systems}
We want to describe the holomorphic solutions of $\mathcal{G}^k$, i.e.
the morphisms of $\mathcal{D}_{\mathbb{C}}^{\an}$-modules
$(\mathcal{G}^k)^{\an}\to\mathcal{O}_{\mathbb{C}}^{\an}$.

Let $\Sigma=\Sigma_1\cup\Sigma_2$. We recall that
$(\mathcal{G}^k)^{\an}(V)$ is isomorphic to $\mathcal{O}^{\an}(V)\otimes_{\mathbb{C}[t]}H^k(\Gamma_{f^{-1}(\mathbb{C}\setminus\Sigma)}(\Omega^{\bullet+n-1}_{U|\mathbb{A}^1}),\nabla)$ for any open subset $V$ in $\mathbb{C}\setminus\Sigma$ (cf. Proposition \ref{prop:analytG}) and $\widetilde{\Hrd_{k+n-1}}(V)$ is the space of families on $V$ of cycles in $\Hrd_{k+n-1}(f^{-1}(t),e^{-g_t})$ which depend continuously on $t$.

\begin{theo}\label{solutions}
Let $V$ be a simply connected open subset of $\mathbb{C}\setminus\Sigma$. The morphism 
$$\begin{array}{llll}
  \Psi:&\widetilde{\Hrd_{k+n-1}}(V)&\to&\Hom_{\mathcal{D}^{\an}(V)}((\mathcal{G}^k)^{\an}(V),\mathcal{O}^{\an}(V))\\
  &(c_t\otimes e^{g_t})_{t\in V}&\mapsto&\big(a\otimes [w]\mapsto(I:t\mapsto\int_{c_t}a(t)w_{|f^{-1}(t)}~e^{g_t})\big)
\end{array}$$
is well-defined and is an isomorphism.

In other words, the space of
holomorphic solutions of $\mathcal{G}^k$ on $V$ is isomorphic to $\widetilde{\Hrd}_{k+n-1}(V)$.
\end{theo}

\begin{proof}

\noindent
\begin{enumerate}
\item As $\nabla_t(a(t)w_{|f^{-1}(t)})=a(t)\nabla(w)$ and as $\int_{\cdot}\cdot~e^{g_t}$ defined a pairing between $\Hrd_{k+n-1}(f^{-1}(t),e^{-g_t})$ and
  $\Hdr^{k+n-1}(f^{-1}(t),e^{g_t})$ (cf. Proposition \ref{prop:perpair}
  for $U_t=f^{-1}(t)$), $I(t)$ is
  well-defined for any $t\in V$.

\item Now we have to prove that $I$ is holomorphic on $V$ and that we
  have constructed a morphism of $\mathcal{D}^{\an}(V)$-modules.

  Let $t_0\in V$. We consider a small disc
  $S\subset\mathbb{C}\setminus\Sigma$ centered at
  $t_0$ such that i) and ii) of Lemma \ref{fibration} are fulfilled.
  We denote by $\widetilde{\phi}$ the restriction of $\phi$ to
  $f^{-1}(t_0)\times S$. It is a homeomorphism from $f^{-1}(t_0)\times
  S$ to $f^{-1}(S)$.

  Let $t\in S$. Then,
$$\begin{array}{lll}
  I(t)-I(t_0)&=&\int_{c_t}a(t)w_{|f^{-1}(t)}e^{g_t}-\int_{c_{t_0}}a(t_0)w_{|f^{-1}(t_0)}e^{g_{t_0}},\\
  &=&\int_{c_t-c_{t_0}}(a\circ f)we^g.
\end{array}$$

Let $C=\underset{s\in[t_0,t]}{\cup}c_s$. Then $\partial C=c_t-c_{t_0}$
and according to the Stokes' Formula, $ I(t)-I(t_0)=\int_Cd((a\circ
f)we^g)$.

But as $\nabla([w])=0$, there exists
$\eta\in\Gamma_{f^{-1}(\mathbb{C}\setminus\Sigma)}(\Omega^{k+n-1}_U)$
such that $dw+dg\wedge w=df\wedge\eta$. Then,
$$\begin{array}{lll}
  d((a\circ f)we^g)&=&\big(df\wedge((\frac{\partial a}{\partial t}\circ f)w)+(a\circ f)(dw+dg\wedge w)\big)e^g,\\
  &=&df\wedge\big((\frac{\partial a}{\partial t}\circ f)w+(a\circ f)\eta\big)e^g,
\end{array}$$
and by Fubini's Theorem, 
$$I(t)-I(t_0)=\int_{t_0}^t\bigg(\int_{c_s}(\big(\frac{\partial a}{\partial t}\circ f)w+(a\circ f)\eta\big)e^g\bigg)ds.$$

Set $\alpha=(\frac{\partial a}{\partial t}\circ f)w+(a\circ f)\eta$.
We have $\int_{c_s}\alpha
e^g=\int_{c_{t_0}}\widetilde{\phi}^{\ast}(\alpha e^{g_t})$.

According to Lemma \ref{fibration} ii), there exists a compact $K$ in $U$ such that $\widetilde{\phi}^{\ast}(\alpha e^{g_t})=\widetilde{\phi}^{\ast}(\alpha)e^{g_{t_0}}$ on $c_{t_0}\setminus K$.
\begin{itemize}
\item As $\widetilde{\phi}^{\ast}(\alpha e^{g_t})$ is continuous in $s$ and $K$ is compact, $\int_{c_{t_0}\cap K}\widetilde{\phi}^{\ast}(\alpha e^{g_t})$ is continuous in $s$.
\item As $e^{g_{t_0}}$ does not depend on $s$ and has rapid decay on $c_{t_0}\setminus K$, the integral $\int_{c_{t_0}\setminus K}\widetilde{\phi}^{\ast}(\alpha e^{g_t})=\int_{c_{t_0}\setminus K}\widetilde{\phi}^{\ast}(\alpha)e^{g_{t_0}}$ is continuous on $s$. 
\end{itemize}

We conclude that $\int_{c_{t_0}}\widetilde{\phi}^{\ast}(\alpha e^{g_t})$ is continuous on $s$ and 
$$\lim_{t\to t_0}\frac{I(t)-I(t_0)}{t-t_0}=\int_{c_{t_0}}\alpha e^g.$$

Then, $I$ is holomorphic on $V$. 

Moreover, as
$\partial_t(a\otimes[w])=\frac{\partial a}{\partial
  t}\otimes[w]+a\otimes[\eta]$, the morphism $\big(a\otimes
[w]\mapsto(I:t\mapsto\int_{c_t}a(t)w_{|f^{-1}(t)}e^{g_t})\big)$ is a
morphism of $\mathcal{D}^{\an}(V)$-modules.

\item At last, we have to prove that $\Psi$ is an isomorphism.
  According to Theorem \ref{thm:dual}, we have an isomorphism:
$$\begin{array}{lll}
  \Hrd_{k+n-1}(f^{-1}(t),e^{-g_t})&\to&\Hom_{\mathbb{C}}(\Hdr^{k+n-1}(f^{-1}(t),e^{g_t}),\mathbb{C})\\
  c_t\otimes e^{g_t}&\mapsto&(w_t\mapsto\int_{c_t}w_te^{g_t}).
\end{array}$$
Then, 
$\widetilde{\Hrd_{k+n-1}}(V)\simeq\Hom_{\mathbb{C}}(\widetilde{\Hdr^{k+n-1}}(V),\mathbb{C})$. 

As
$\Hom_{\mathbb{C}}(\widetilde{\Hdr^{k+n-1}}(V),\mathbb{C})\simeq\Hom_{\mathcal{D}^{\an}(V)}((\mathcal{G}^k)^{\an}(V),\mathcal{O}^{\an}(V))$
(cf. Corollary 7.1.1 p. 71 of \cite{pham}), we conclude that $\Psi$ is
an isomorphism.
\end{enumerate}
\end{proof}

\end{document}